\documentclass[a4paper,12pt]{amsart}

\usepackage{preprint}
\usepackage{amsthm}
\usepackage{amsfonts}

\usepackage{amsmath}

\usepackage[english]{babel}
\usepackage{fancybox}

\nonstopmode

\title{}

\author{}

\usepackage{psfig}

\newcounter{minutes}
\divide\time by 60
\newcounter{hours}
\multiply\time by 60
\addtocounter{minutes}{-\time}


\begin{document}

\begin{center}
{\Large \bf   Experiments with moduli of quadrilaterals}
\end{center}
\medskip

\begin{center}
{\large \bf  Antti Rasila and Matti Vuorinen}
\end{center}

\bigskip

\begin{center}
Published in Rev. Roumaine Math. Pures Appl.  51 (2006), 
5--6, 747--757.
\end{center}

\subsection*{Abstract}
Basic facts and definitions of conformal moduli of rings
and quadrilaterals are recalled.
Some computational methods are reviewed.
For the case of quadrilaterals with polygonal sides,
some recent results are given. Some numerical experiments
are presented. This paper is based on 
\cite{BSV} and \cite{RV}.

\section{Introduction}

We give a brief introduction to the conformal moduli of
quadrilaterals and rings. For a comprehensive survey of this
topic see \cite{Kuhnau:2005}.

The capacity of condensers has been studied because of
its importance in physics and its close relation with the potential 
theory  and the  theory of conformal and quasiconformal mappings. The 
analytic computation of capacity is possible only for very few types of 
condensers and for this reason several methods have been developed for 
the numerical computation of capacity. 

Let $E$ and $F$ be two disjoint compact sets in the extended complex
plane ${\CC}$. We assume that each of $E$ and $F$ is the union
of a finite number of nondegenerate disjoint continua, and that the open
set $R={\CC}\setminus (E\cup F)$ is connected. Without loss 
of generality, we also assume that $\infty\notin E$. The domain $R$ is
a {\it condenser}. The complementary compact sets $E$ and $F$ are the
{\it plates} of the condenser. The {\it capacity} of $R$ is defined by
\begin{equation}
\label{cap}
{\rm cap}\,R=\inf_u \int_R |\triangledown u|^2\,dm,
\end{equation}
where the infimum is taken over all nonnegative, piecewise
differentiable functions $u$
with compact support in $R\cup E$ such that $u=1$ on $E$.
It is well-known that under the assumptions we made above, $R$ is 
regular
for the Dirichlet problem and the harmonic function on $R$ with
boundary values $1$ on $E$ and $0$ on $F$ is the unique function that
minimizes the integral in (\ref{cap}). This function is called the {\it
potential function} of the condenser.

Capacity is a conformal invariant: Suppose that $f$ maps $R$ conformally
onto $R^\prime$. Let $E$ and $F$ correspond to $E^\prime$ and $F^\prime$
respectively (in the sense of the boundary correspondence under 
conformal
mapping). Then cap$\,R={\rm cap}\,R^\prime$. This property can be used 
for
the
analytic computation of capacity provided that the capacity of some
`canonical' condensers is known and the corresponding conformal mappings
can be constructed. Unfortunately such an analytic computation can be 
made
only for very few doubly-connected condensers; see \cite{IvTr}.

If both $E$ and $F$ are connected (and hence $R$ is doubly-connected), 
$R$ is called a {\it ring domain}. A ring domain $R$ can be mapped 
conformally onto the annulus 
$\{z:1<|z|<e^M\}$, where $M=\mathrm{mod}\,R$ is the {\it conformal 
modulus} 
of the ring 
domain $R$, defined by $\mathrm{mod}\,R=2\pi/\mathrm{cap}\,R$.
See also \cite{Ah}, \cite{Hen}, \cite{Kuhnau:2005}.

A Jordan domain $D$ in $\C$ with marked points 
$z_1,z_2,z_3,z_4\in \partial D$ is a quadrilateral and denoted by 
$(D;z_1,z_2,z_3,z_4)\, .$ We use the canonical map 
onto a rectangle $(D';0,1,1+ih,ih)$ to define the modulus $h$ of a 
quadrilateral $(D;z_1,z_2,z_3,z_4)\,.$ The modulus of  
$(D;z_2,z_3,z_4,z_1)\,$ is $1/h \, .$ We mainly study the situation
where the boundary of $D$ consists of the polygonal line segments
through $z_1,z_2,z_3,z_4$ (always positively oriented). In this case, the 
modulus is denoted by $\mathrm{QM}(D;z_1,z_2,z_3,z_4)$. If the boundary
of $D$ consists of straight lines connecting the given boundary points,
we omit the domain $D$ and denote the quadrilateral and the
corresponding modulus simply by $(z_1,z_2,z_3,z_4)$
and $\mathrm{QM}(z_1,z_2,z_3,z_4)$. 

\medskip

\centerline{
\psfig{figure=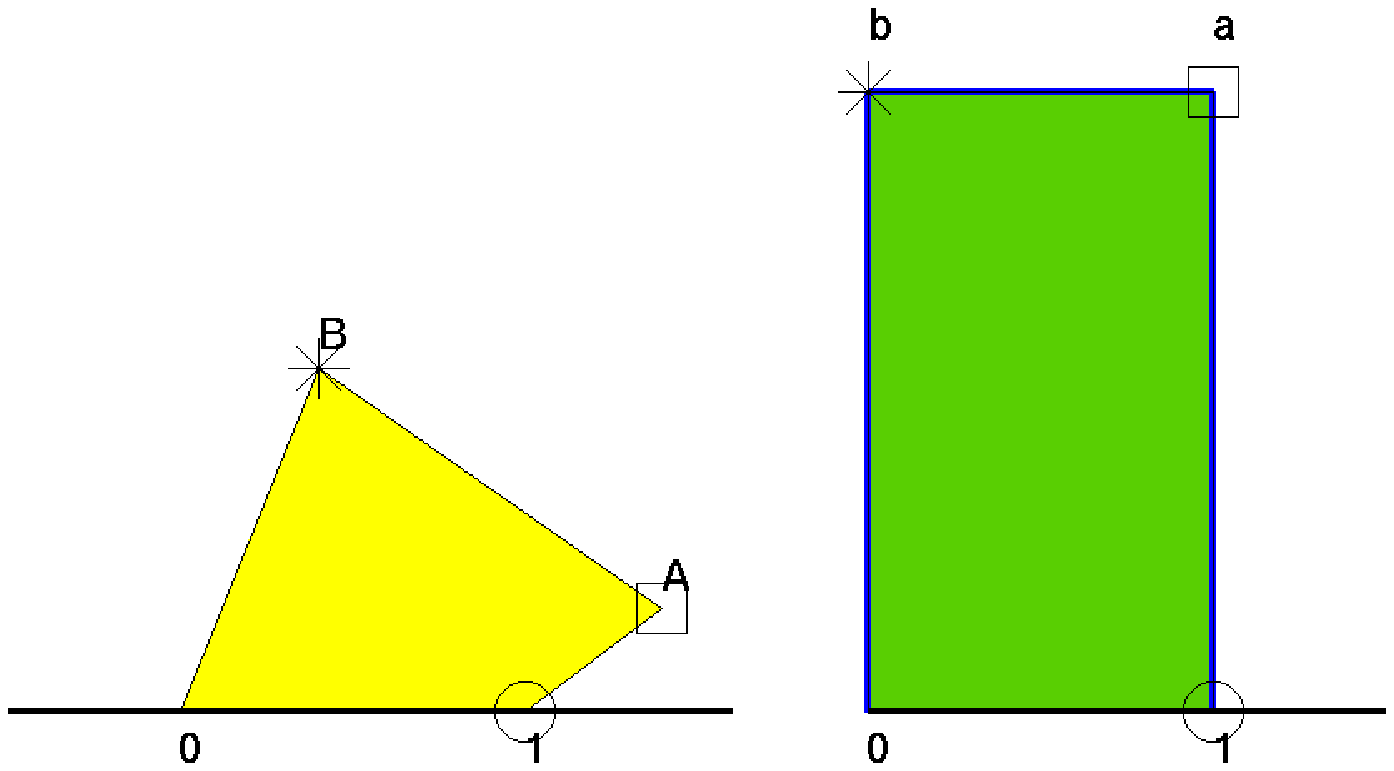,width=11cm}
}
\begin{capti}
The height of the canonical rectangle is $\mathrm{QM}(A,B,0,1).$
\end{capti}

\medskip

The following problem is known as the {\it Dirichlet-Neumann problem}.
Let $D$ be a region in the complex plane whose boundary
$\partial D$ consists of a finite number of regular Jordan
curves, so that at every point of the boundary a normal is
defined. Let $\psi$ to be a real-valued continuous
function defined on $\partial D$. Let $\partial D =A \cup B$ where
$A, B$ both are unions of Jordan arcs. Find a function $u$ satisfying
the following conditions:
\begin{enumerate}
\item
$u$ is continuous and differentiable in
$\overline{D}$.
\item
$u(t) = \psi(t),\qquad t \in A$
\item
If $\partial/\partial n$ denotes differentiation in
the direction of the exterior normal, then
$$
\frac{\partial}{\partial n} u(t)=\psi(t),\qquad t\in  B.
$$
\end{enumerate}

 One can express the modulus of a quadrilateral
$(D; z_1, z_2, z_3, z_4)$ in terms of the solution of the
Dirichlet-Neumann problem as follows. Let $\gamma_j, j=1,2,3,4$ be the 
arcs of
$\partial D$ between $(z_1, z_2)\,,$ $(z_2, z_3)\,,$ $(z_3, z_4)\,,$
$(z_4, z_1),$ respectively. If $u$ is the (unique) harmonic solution of 
the 
Dirichlet-Neumann problem with boundary values equal to 0 on $\gamma_2$, 
equal to 1 on $\gamma_4$ and with $\partial u/\partial n = 0$ on 
$\gamma_1 \cup \gamma_3\,,$ then by \cite[p.65/Thm 4.5]{Ah}:
\begin{equation}
\label{qmod}
\mathrm{QM}(D;z_1,z_2,z_3,z_4)=
\int_D |\triangledown 
u|^2\,dm.
\end{equation}

In conclusion, the computation of the modulus of a ring or 
quadrilateral can be reduced to solving the  Dirichlet problem 
(\ref{cap}) or the Dirichlet-Neumann problem, respectively.
The connection between ring and quadrilateral moduli is given
in \cite[p.102]{Kuhnau:2005} or \cite[p.36]{Lehto:1973}.

\section{Review of some numerical methods}

Both for the  ring  and the quadrilateral
case we may consider the following methods:
\begin{itemize}
\item[(a)] Approximate the canonical conformal map.
\item[(b)] Approximate the solution of the Dirichlet 
(or Dirichlet-Neumann) problem.
\end{itemize}

The recent survey of Wegmann \cite{Weg05} provides an extensive
review of the various techniques of the approximation of conformal
mappings. See also Driscoll and Trefethen \cite{DrTr}, Papamichael
\cite{Pap}. We now mention some of the known methods, following
closely \cite{BSV}.

The paper \cite{Gai2} of D.~Gaier includes a review of the various 
methods applicable to the computation of the capacity of planar ring 
domains.

The finite element method was first applied to the computation of
capacity by G.~Opfer \cite{Opf}. Several numerical experiments are
reported by J.~Weisel \cite{Wei2}. Another numerical method
is based on the Gauss-Thompson principle which implies a formula for the
capacity involving Green's function. Numerical computations are given in
\cite{Wei1}. 
N.~Papamichael and his collaborators \cite{PaKo}, \cite{PS}, 
\cite{PaWa} have
developed an orthonormalization technique for the approximation of the
conformal mapping of doubly-connected domains. This technique gives, in
particular, approximations of capacity. Many numerical
computations are presented in the above papers. 

The capacity of a polygonal ring domain can be also computed by the
Schwarz-Christoffel transformation which provides a semi-explicit
formula for the conformal mapping of the domain onto an annulus
(see \cite{Hen}). For simply-connected domains this methods has 
been developed by T. Driscoll, L.N.~Trefethen and their coauthors 
(see  \cite{DrTr}, \cite{DrVa}). 
It seems that for doubly-connected polygonal domains
the only related works are those of H.~Daeppen \cite{Dae} and C.~Hu
\cite{Hu}. Hu's method has been tested successfully in several 
computations, (see \cite{Hu}, \cite{BeVu}). It is partially based
on the wise choice of certain points on the complementary sets $E$, $F$
of the ring domain.

Another numerical-analytic method that can be used for the computation
of capacity is the multipole method. The potential function is written
as a linear combination of explicit basic functions (multipoles) with
unknown coefficients. The coefficients are then computed numerically.
The multipoles constitute a complete, minimal system in a certain
Hardy-type space of functions. Their construction is based on the theory
of conformal mapping. This method has been developed by V.I.~Vlasov (see 
\cite{Vla} 
and references therein) as a general method for numerical solution of a 
wide class 
of boundary value problems. He has applied this approach to find the 
potential function of condensers.

\section{Web-based simulations}

The solutions of the Dirichlet and the Dirichlet-Neumann problems can 
be 
approximated by the method of finite elements, see \cite[pp. 
305--314]{Hen}, \cite{Pap}. Hence, this method can also be 
used to approximate the modulus of quadrilaterals and rings.
The Dirichlet-Neumann problem can be numerically solved with AFEM
(Adaptive FEM) numerical PDE analysis package by Klas Samuelsson. This software 
applies, e.g., to multiply connected polygonal domains. In particular,
we may use it to compute the modulus (capacity) of a bounded ring whose
boundary components are broken lines. Examples and applications
for  this software are given in \cite{BSV}. In \cite{HVV} a theoretical 
formula for computing $\mathrm{QM}(A,B,0,1)$ was given with its
implementation with Mathematica. This lead to a
more systematic study of the modulus of quadrilateral
in \cite{DV}. In the course of the work on \cite{DV} several
conjectures were formulated and this lead us to
look for an improved version of the algorithm in \cite{HVV}
for the computation of $\mathrm{QM}(A,B,0,1)$.
It seems that the AFEM software of Samuelsson is very efficient for this 
purpose.

Our goal is to use the AFEM software of Samuelsson for computations
involving moduli of polygonal rings and quadrilaterals, 
and to write a user-interface
providing access to AFEM via a web browser (e.g. Mozilla).
The advantages in this 
approach are: 
(1) no programming needed to use AFEM and (2) mobile computing: 
available for  
everyone. Currently the pilot tests work locally at the local university 
network.

\medskip

\centerline{
\psfig{figure=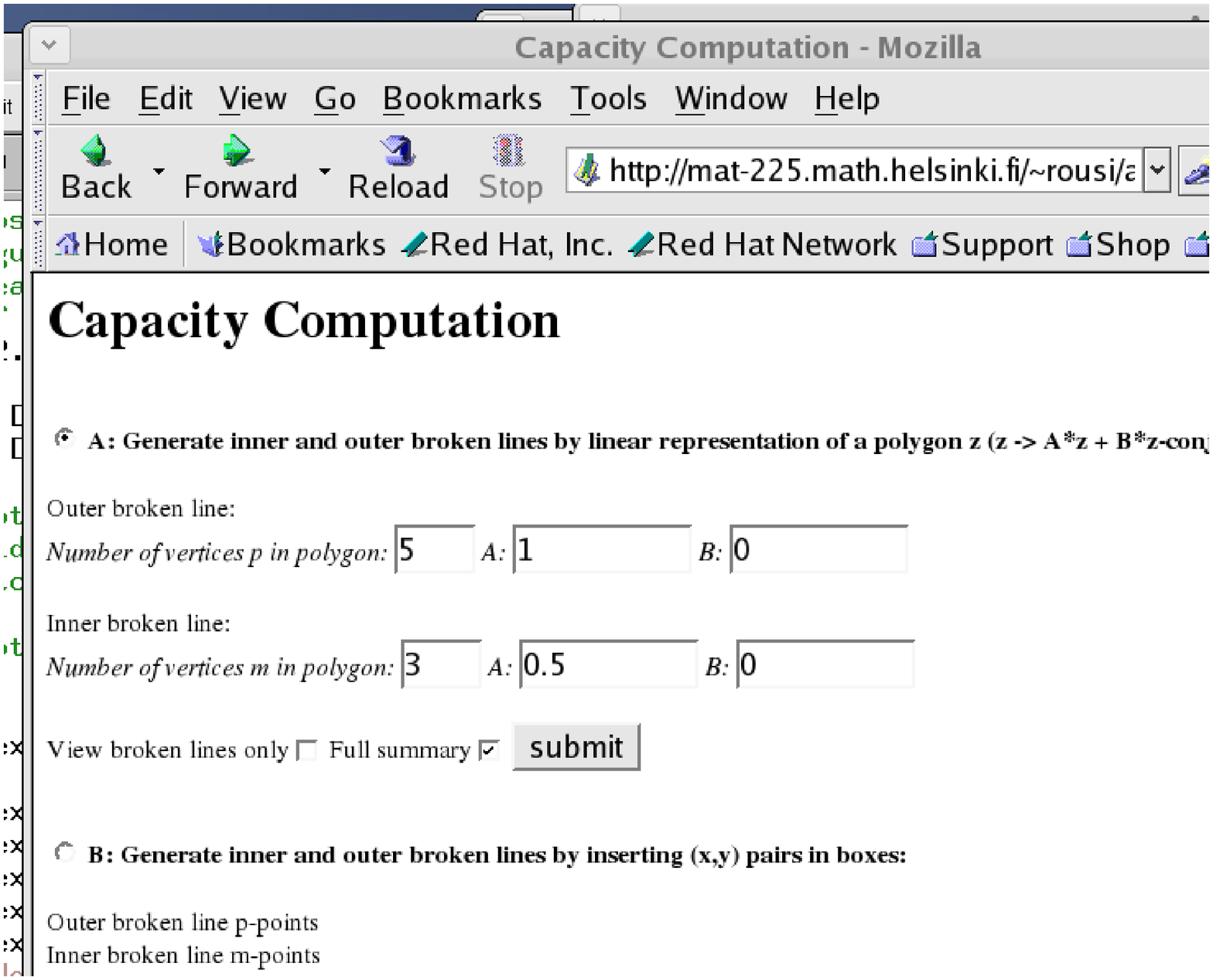,width=12cm}
}
\begin{capti}
Entering two regular polygons for computing the capacity of 
the corresponding ring domain by using a web browser.
\end{capti}

\medskip

\centerline{
\psfig{figure=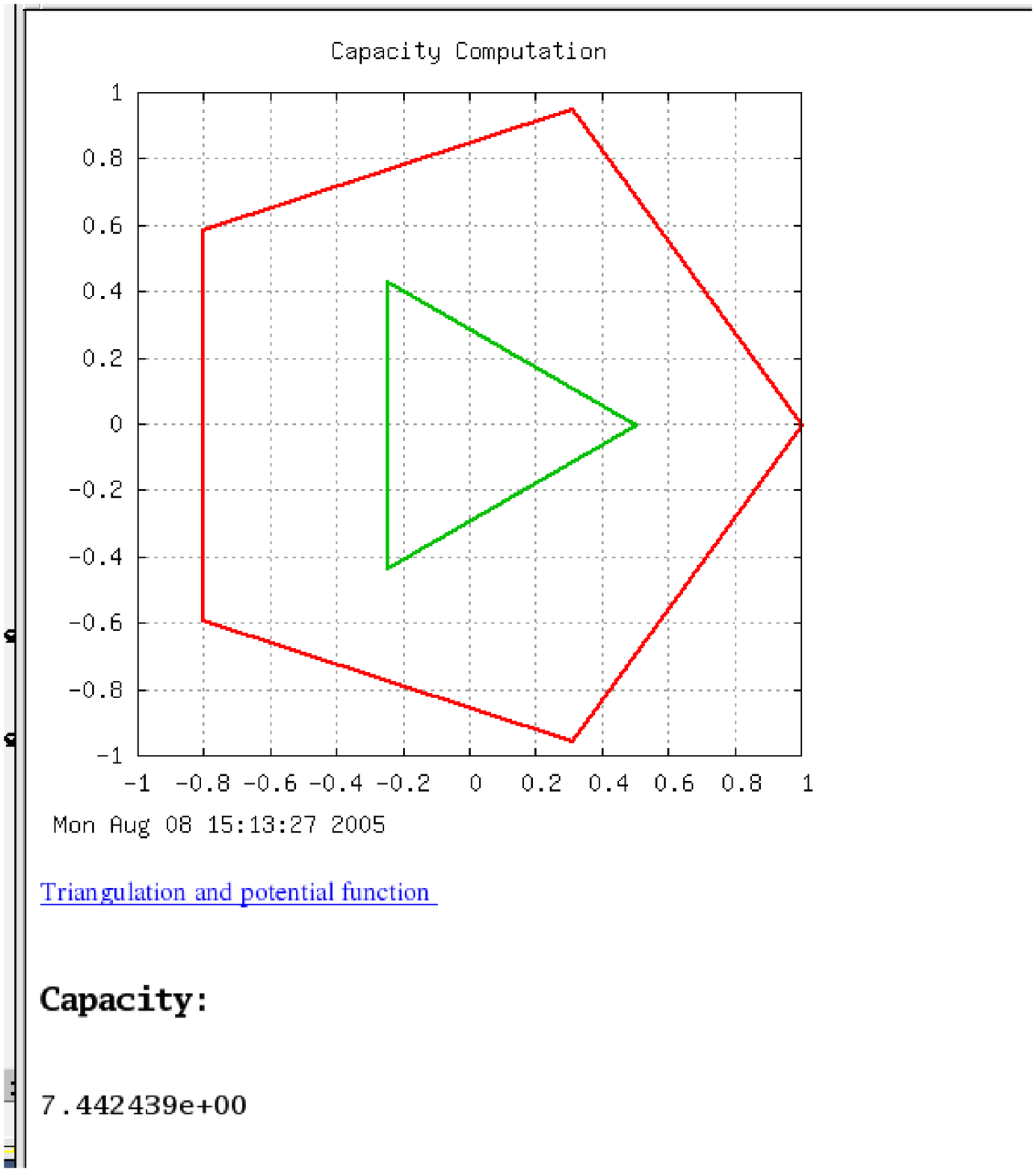,width=8cm}
}
\begin{capti}
The output of the program is also displayed by using the web-based
interface.
\end{capti}

\section{Experiments}

We give some examples of problems which can be studied by using the
AFEM software.

\begin{exmp}
Let $a,b \in C$ with $\im a >0, \im b >0$ and assume that $ (a,b,0,1)$
determines the vertices of a quadrilateral and
$ \arg b \in (\pi/2,\pi),$ $ \arg (a-1) \in (0,\pi/2).$
Is it true that
\begin{equation}
\label{trans}
 \mathrm{QM}(a,b,0,1) \le \mathrm{QM}(1+i|a-1|,i|b|,0,1) ?
\end{equation}
\end{exmp}

\medskip

\centerline{
\psfig{figure=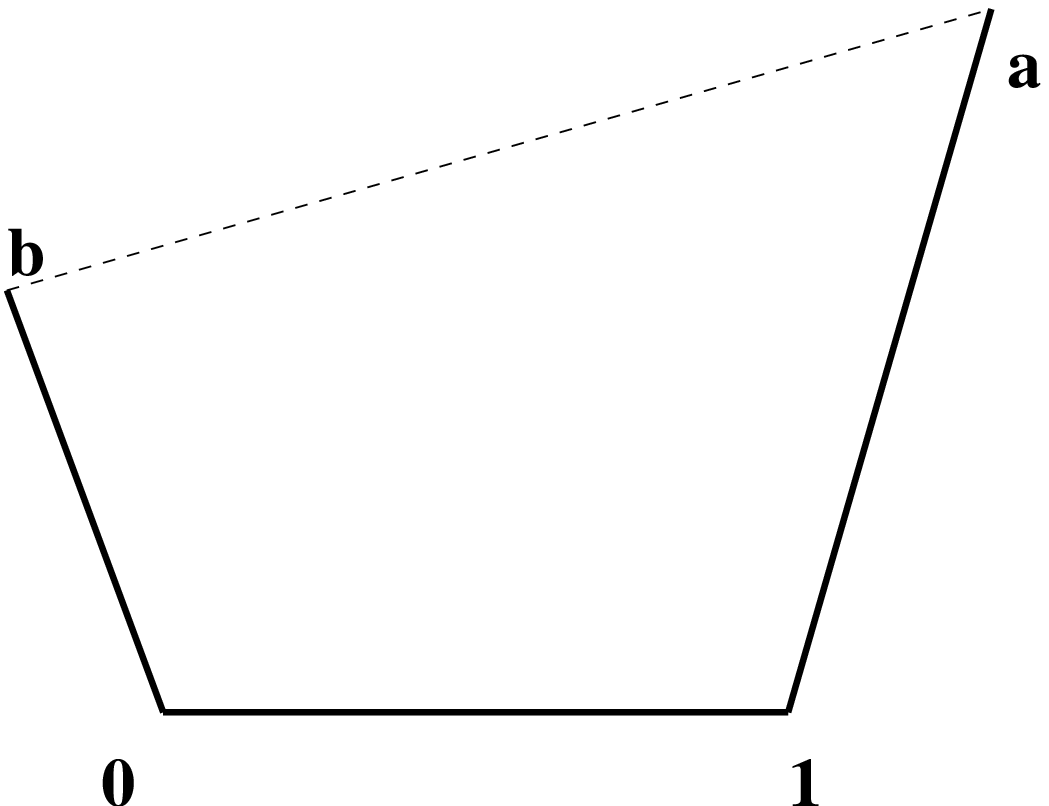,width=5.5cm}
\qquad\qquad 
\psfig{figure=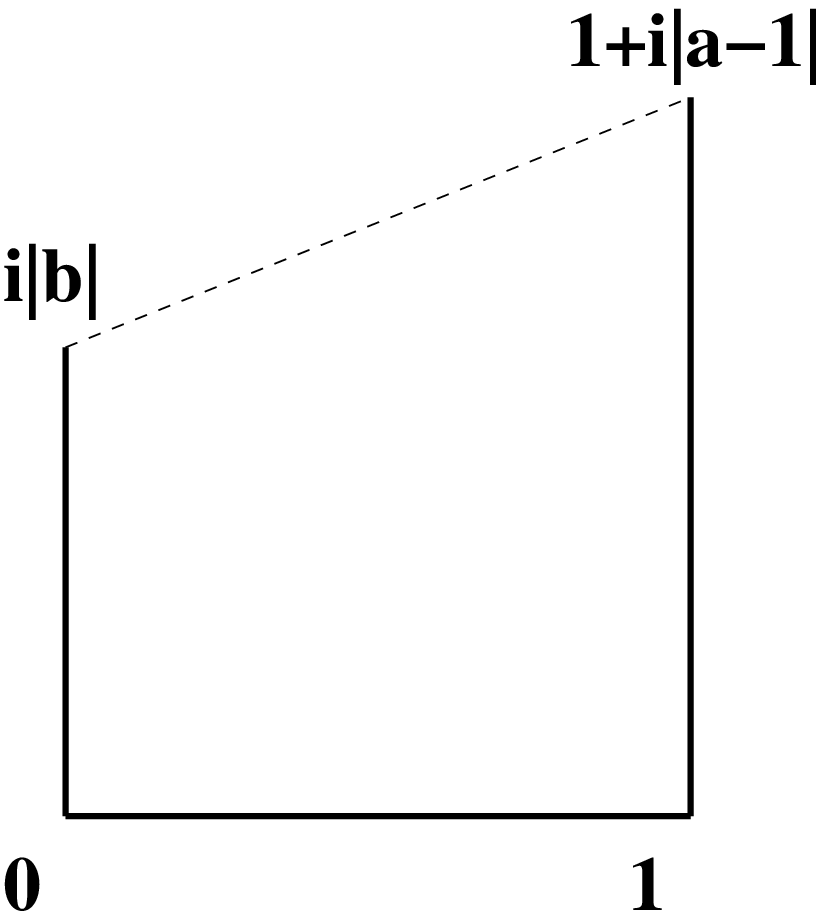,width=4.5cm}
}
\begin{capti}
The quadrilaterals $(a,b, 0,1)$
and $(1+|a-1|i, |b|i,0,1)$ for the left and
right sides for (\ref{trans}), respectively.
\end{capti}

\medskip

\centerline{
\psfig{figure=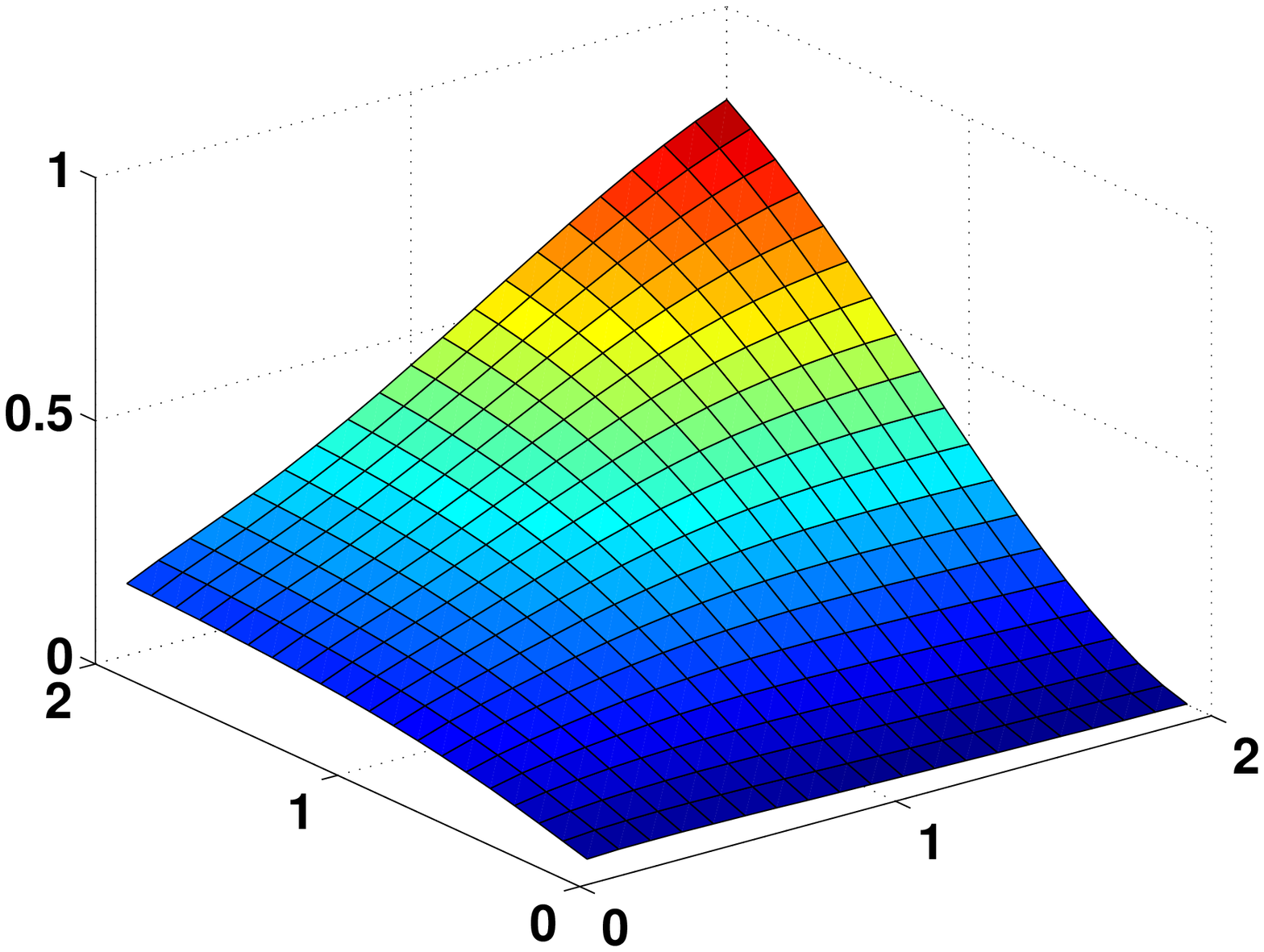,width=7cm}
\qquad
\psfig{figure=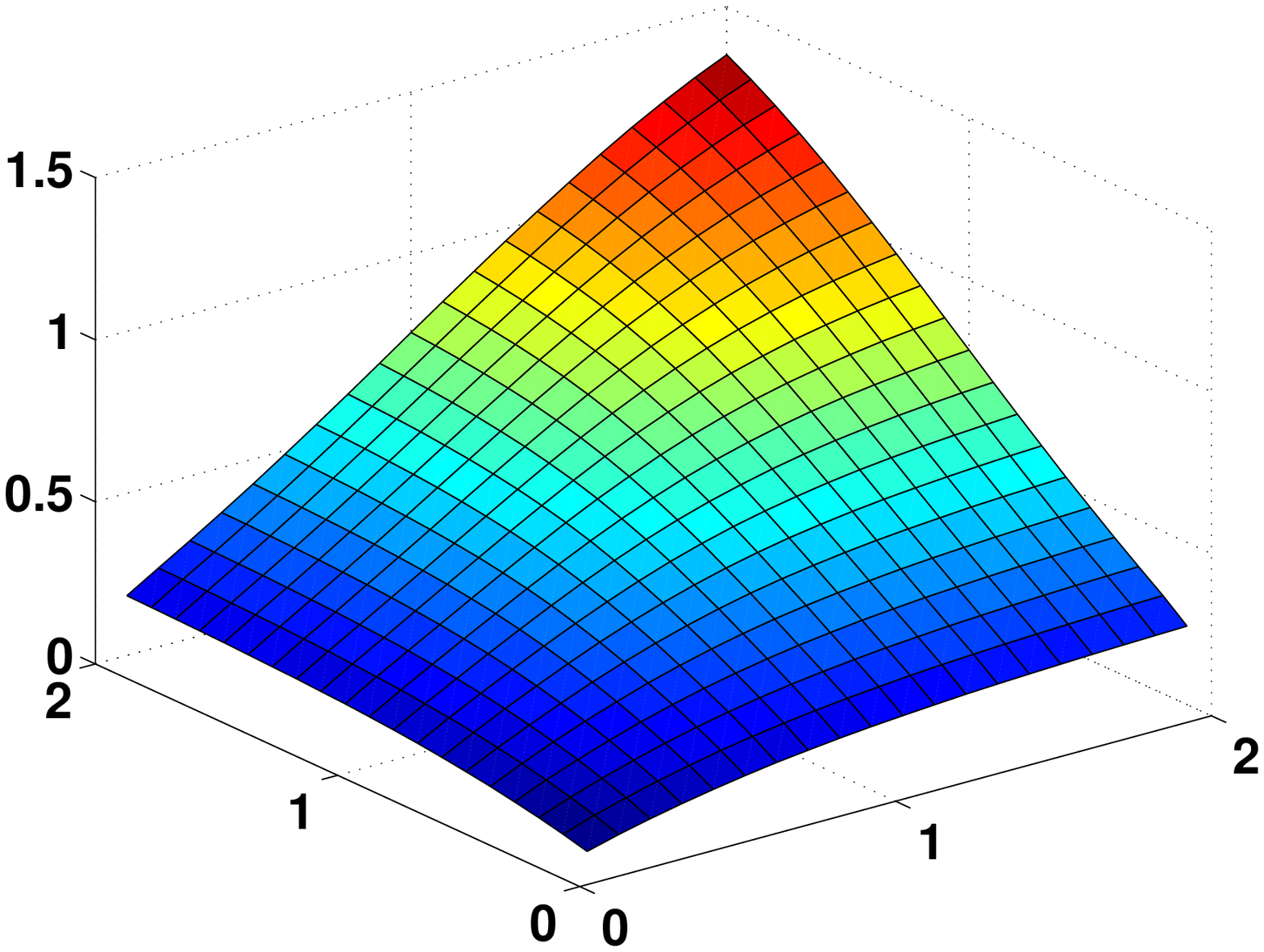,width=7cm}
}

\begin{capti} 
The function 
$f(x,y)=\mathrm{QM}(1+xi,yi,0,1)-\mathrm{QM}(1+x\exp(i\alpha),y\exp(i\beta),0,1)$
for $x,y\in (0,2)$. Here $\beta=3\pi/4$ and $\alpha$ is $\pi/2$ (left)
and $\pi/8$ (right). 
\end{capti}

\bigskip

\begin{remark}
\label{mrem}
The expression on the right hand side of (\ref{trans}) has an 
analytic
expression if $|a-1|= h= |b|+1\,.$
Bowman \cite[pp.\ 103-104]{Bo} gives a formula for the conformal
modulus of the quadrilateral with vertices  $1+hi$,  $(h-1)i$, $0$, and 
$1$
when $h>1$ as $M(h)\equiv\K(r)/\K(r')$ where
$$
r = \bigg(\frac{t_1-t_2}{t_1+t_2}\bigg)^2\, , \quad
t_1 = \mu^{-1}\left(\frac{\pi}{2c}\right) \, , \quad
t_2 = \mu^{-1}\left(\frac{\pi c}{2}\right)\, , \quad
c = 2h - 1 \, .
$$
Here for $0<r<1$
$$
\mu(r) = \frac{\pi}{2}\frac{\K(r')}{\K(r)} \, , \quad
\K(r)=\int_0^1 \frac{dx}{\sqrt{(1-x^2)(1-r^2 x^2)}}\, ,
$$
and $r'= \sqrt{1-r^2}\,.$

Therefore, the quadrilateral can be conformally mapped onto the
rectangle  $1+iM(h)$, $i M(h)$, $0$, $1$, with the vertices 
corresponding
to each other. It is clear that $h-1 \le M(h) \le h \,.$
The formula $$M(h)= h + c +O(e^{-\pi h}),  c=-1/2- \log 2/\pi\,,$$
is given in \cite{PS}. As fas as we know there is neither
an explicit nor asymptotic formula
for the case when the angle $\pi/4$ of the trapezoid is equal
to $\alpha \in (0, \pi/2) \,.$
\end{remark}

\begin{exmp} (Duplication formula)
Let $\phi \in (0,\pi)$, $h, k >0$, $A= 1+ h \exp(i
\phi)$, $B= k\exp(i \phi)$. We study when the 
following inequality holds:
\begin{multline}
\label{dupl}
\mathrm{QM}(A,B,0,1)+ 
\mathrm{QM}(\overline{(1-B)},\overline{(1-A)},0,1)
\le
\mathrm{QM}(A,B,1-A, 1-B).
\end{multline}
Here equality holds if $A=1+ih$, $B=ih$ and $h>0 \,.$
In this special case the result may be regarded as a 
duplication formula.
\end{exmp}

\medskip

\medskip
\centerline{
\psfig{figure=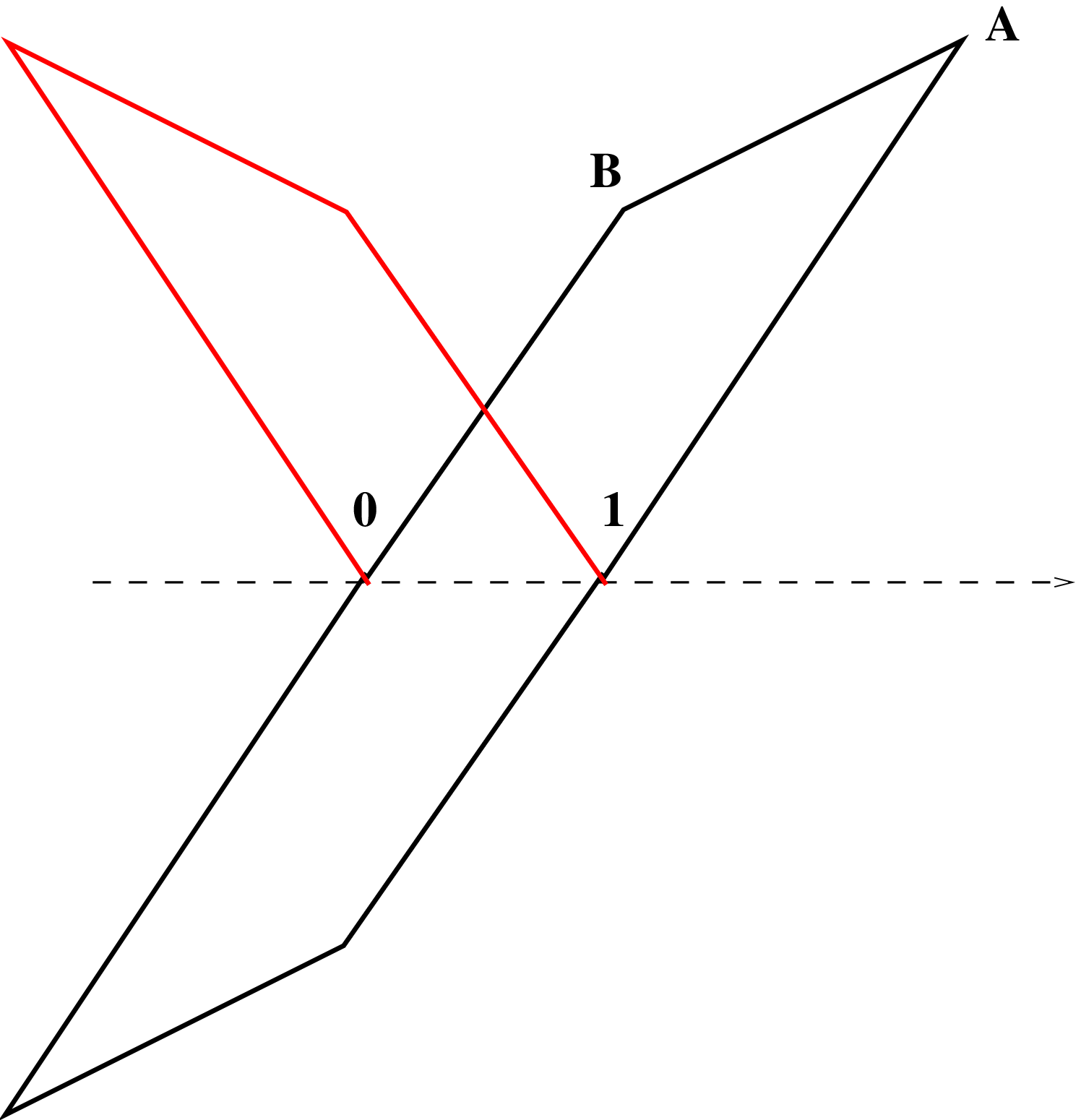,width=5.5cm}
}

\begin{capti}
The quadrilaterals in (\ref{dupl}).
\end{capti}


\medskip

\centerline{
\psfig{figure=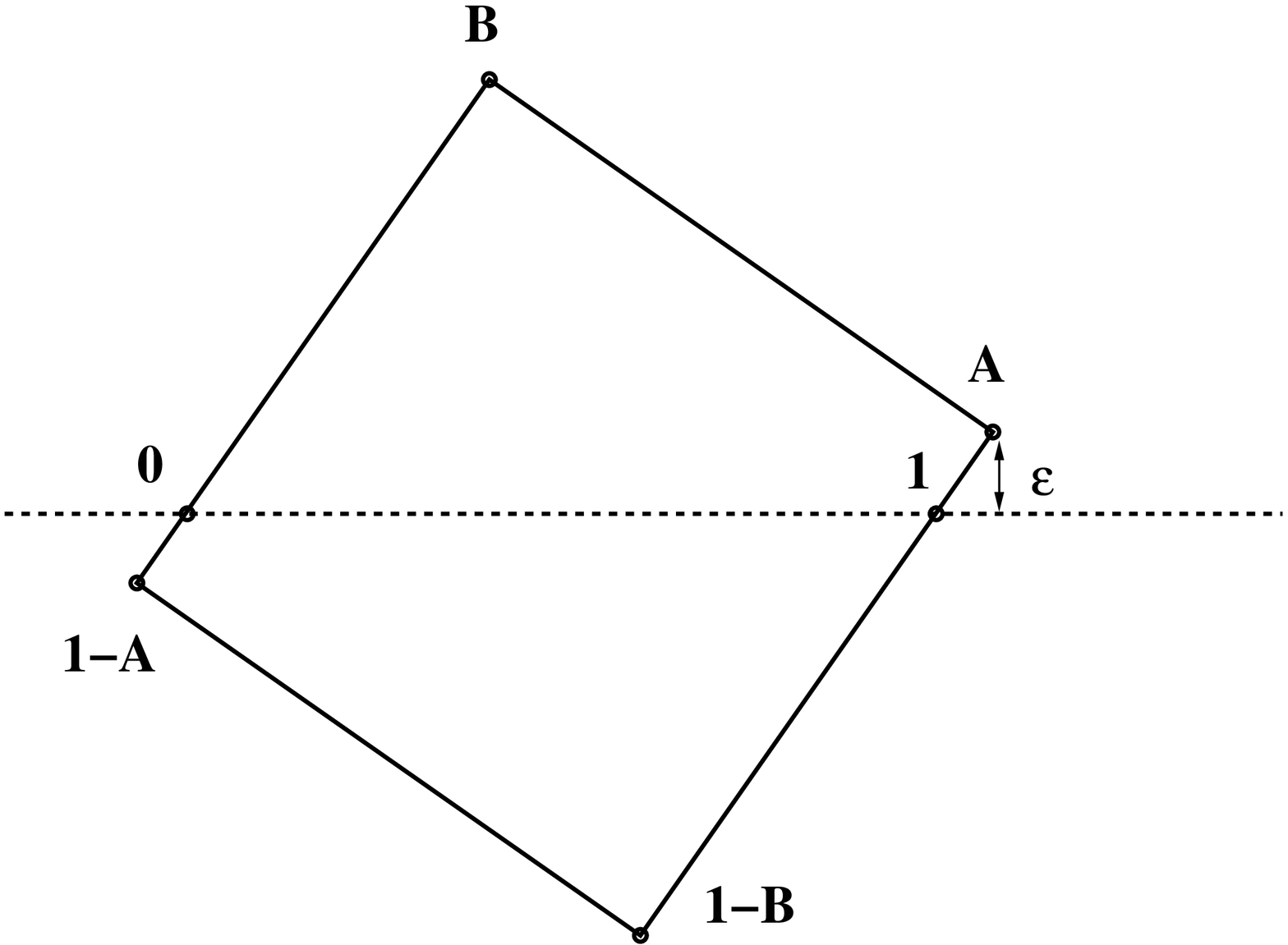,width=7.5cm}
}

\begin{capti}
If $B = 1/2+ i(A-1/2)$, then $A,B, 1-A, 1-B$ are the vertices of
a square and hence $\mathrm{QM}(A,B, 1-A, 1-B)= 1$. Letting 
$|A-1| \to 0$ we see that the left side of (\ref{dupl}) tends to $0$ 
whereas the right side $= 1$.
\end{capti}

\medskip

\centerline{
\psfig{figure=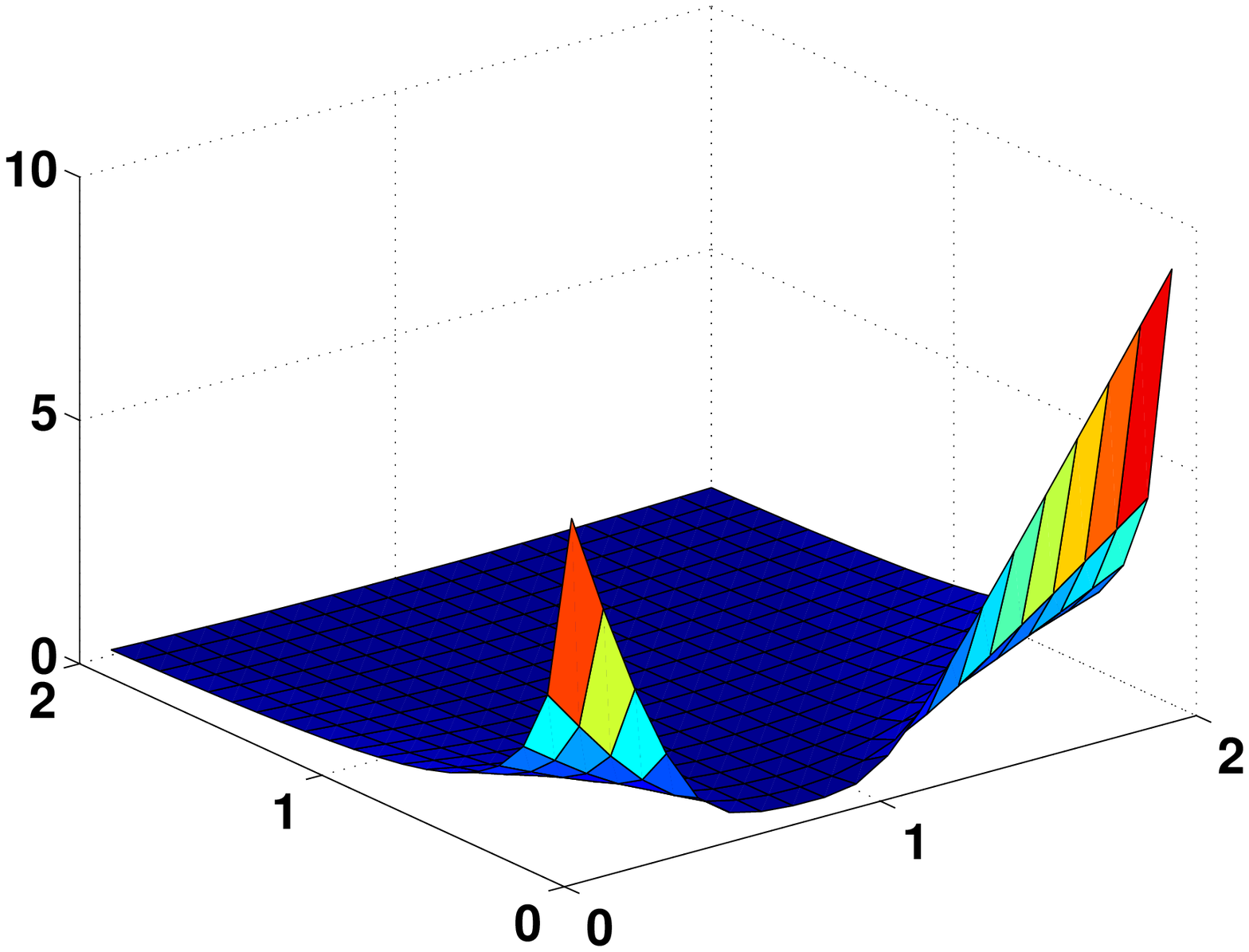,width=7.2cm}
}
\begin{capti} 
The function 
$$g(x,y)=
\mathrm{QM}(A,B,1-A, 1-B)-
\mathrm{QM}(A,B,0,1)-\mathrm{QM}(\overline{(1-B)},\overline{(1-A)},0,1)
$$
where $A=x+iy$ and $B=\exp(i\arg(A-1))$ for $x,y\in (0,2)$.
\end{capti}

\begin{remark}
Let $h,k >1$ and consider next the rectangle with the 
vertices $A=1+i(h+k-1),  B=i(h+k-1), 0, 1 \,.$ If we split this
rectangle in two trapezoids with the segment joining $i(h-1)$ and
$1+ ih$ and apply (\ref{dupl}), then we get
$$ \mathrm{QM}(A,B, 0,1) \ge \mathrm{QM}(A,B,i(h-1), 1+ih)+ 
\mathrm{QM}(0,1,1+ih,i(h-1)) 
\,.$$
If we use the notation from \ref{mrem} and use the
formula $\mathrm{QM}(A,B, 0,1) =h+k-1\,$ then we can express this as
$$ h+k-1 \ge M(h)+M(k) \ge h+k-2 \,.$$
\end{remark}

\begin{exmp}
\label{dvprob}
(Open problem \cite{DV})
Fix $r,s >0, \alpha \in (0, \pi/2), \beta\in(\pi/2, \pi) \,.$
Determine $t>0$ by the condition that the quadrilaterals
$Q_1=(1+2r e^{i \alpha} , 2s e^{i \beta}, 0,1    )$
and $Q_2=(t+ir, is, -is, t-ir)$ have equal areas. Is it true that
\begin{equation}
QM(1+2r e^{i \alpha} , 2s e^{i \beta}, 0,1    ) \le
QM( t+ir, is, -is, t-ir)\, ? 
\end{equation}
\end{exmp}

\medskip

\centerline{
\psfig{figure=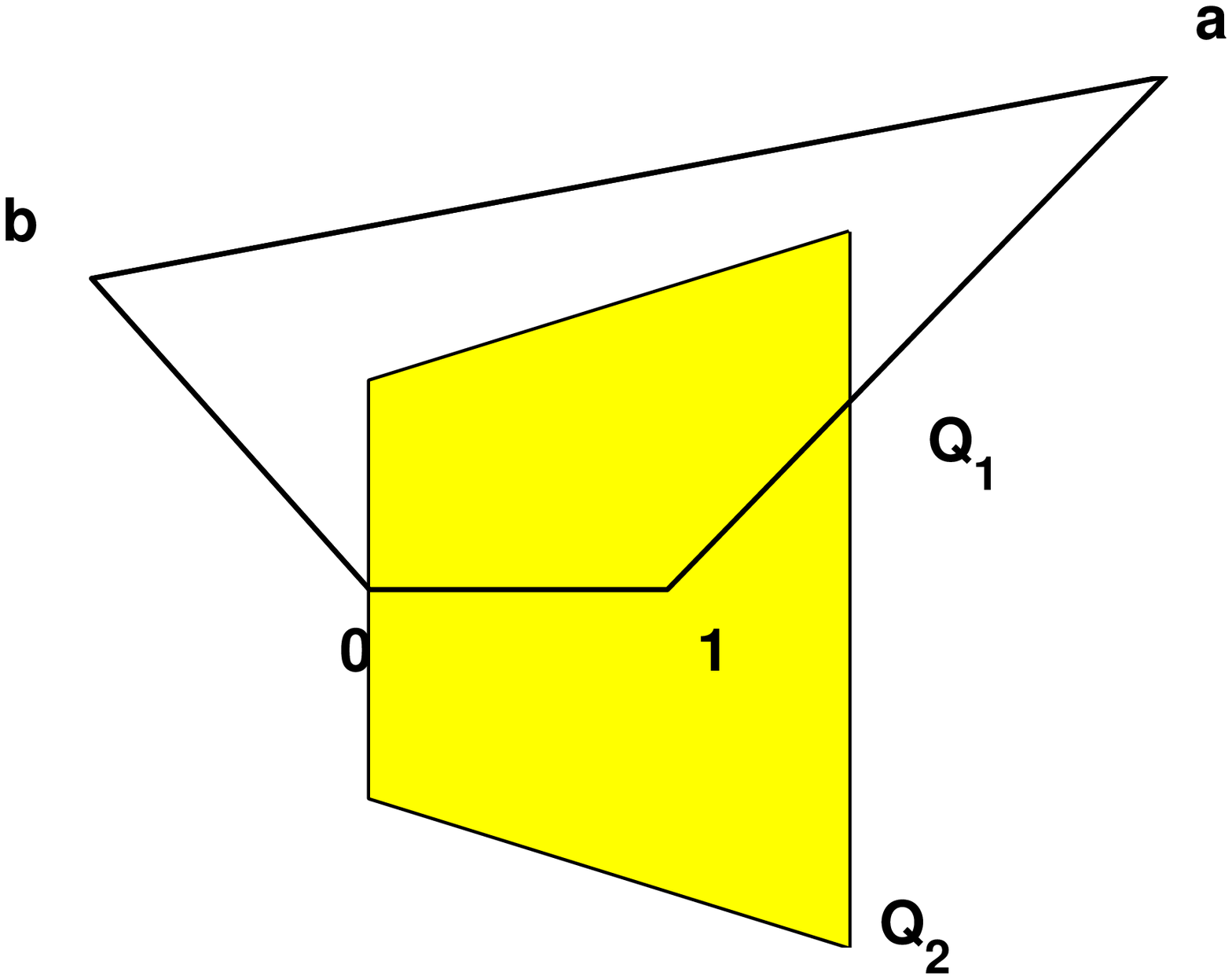,width=5.5cm}
}

\begin{capti}
Quadrilaterals $Q_1$ and $Q_2$ of Example \ref{dvprob}.
\end{capti}

\medskip

\centerline{
\psfig{figure=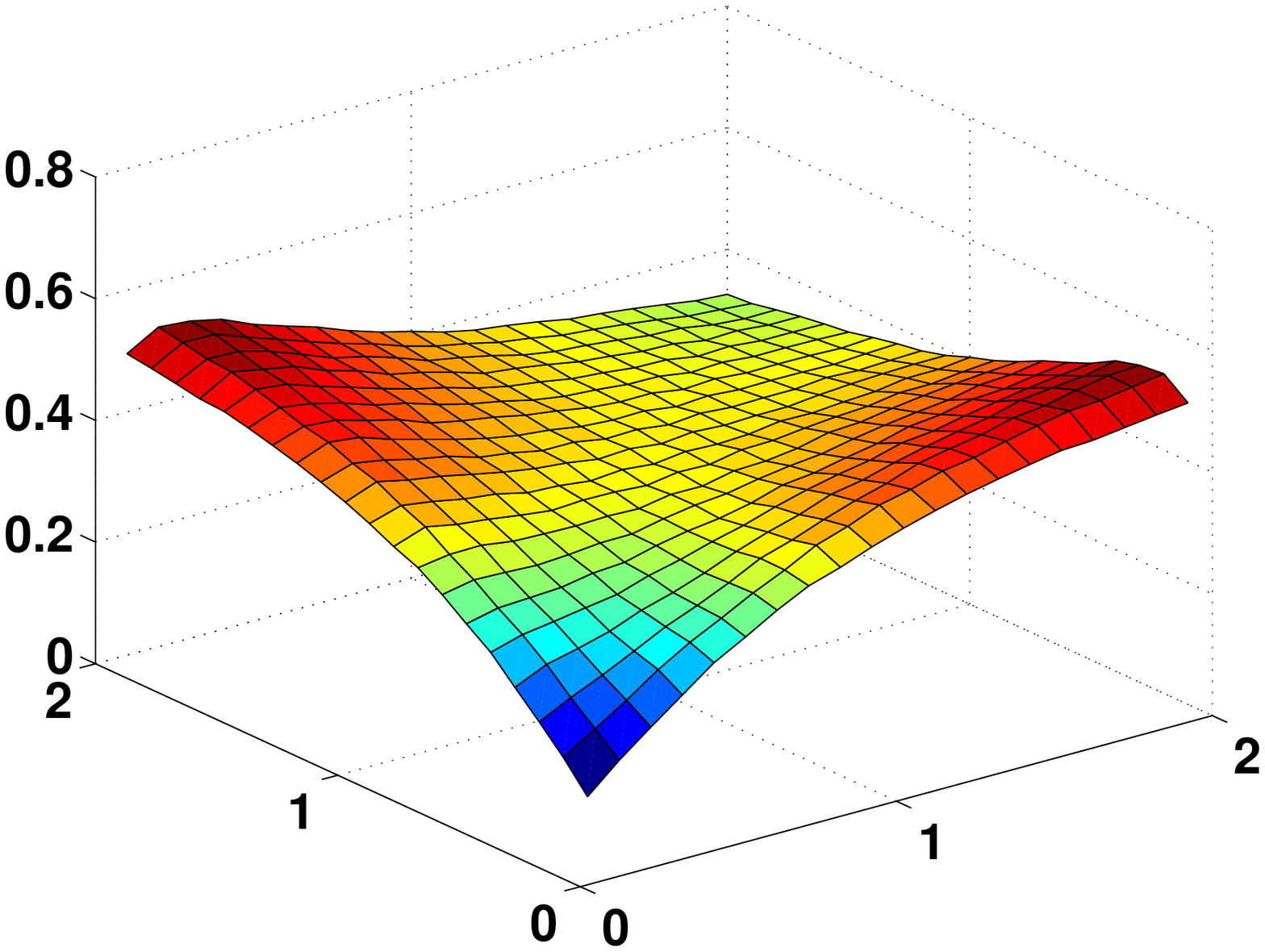,width=7.5cm}
}

\begin{capti}
Function 
$$
h(r,s)=QM( t+ir, is, -is, t-ir) -QM(1+2r e^{i \alpha} , 2s e^{i 
\beta}, 0,1    )
$$
for $\alpha=\pi/4$, $\beta=3\pi/4$.
\end{capti}

\medskip

\subsection*{Acknowledgements}
The research was supported by the Academy of Finland,
by the Romanian Academy and by the Foundation Vilho, Yrj\"o ja Kalle 
V\"ais\"al\"an rahasto of Finnish Academy of Science and Letters.

\bigskip

\noindent
{\bf Antti Rasila:}\\
Helsinki University of Technology\\
Institute of Mathematics\\
P.O.Box 1100, FI-02015 TKK\\
FINLAND \\
Email: {\tt antti.rasila@tkk.fi}

\medskip
\noindent
{\bf Matti Vuorinen:}\\
Department of Mathematics\\
FIN-20014 University of Turku\\
FINLAND\\
Email: {\tt vuorinen@utu.fi}


\begin{thebibliography}{11111}

\bibitem[Ahl]{Ah}
{\sc L. V. Ahlfors,\,}
{\sl Conformal invariants: topics in geometric function theory. }
 McGraw-Hill Book Co.,  1973.


\bibitem[BV]{BeVu} {\sc D. Betsakos} and {\sc M. Vuorinen,\,}
{\sl Estimates of conformal capacity}.
Constr. Approx. 16.4 (2000), 589-602.

\bibitem[BSV]{BSV}
{\sc D. Betsakos, K. Samuelsson} and {\sc M. Vuorinen,\,}
{\sl Computation of capacity of planar condensers},
Publ. Inst. Math. 75 (89) (2004), 233-252.

\bibitem[Bow]{Bo} {\sc F.\ Bowman}: {\sl Introduction to Elliptic
Functions with applications,} English Universities Press Ltd., London,
1953.

\bibitem[Dae]{Dae} {\sc H. Daeppen,\,}
{\sl Die Schwarz-Christoffel-Abbildung f\"ur zweifach zusammenh\"angende
Gebiete mit Anwendungen}. Ph.D. thesis, E.T.H., Zurich, 1988.

\bibitem[DrTr]{DrTr} 
{\sc T.A. Driscoll} and  {\sc L.N. Trefethen:}
{\sl Schwarz-Christoffel mapping.} Cambridge Monographs on Applied and 
Computational Mathematics, 8. Cambridge University Press, Cambridge, 
2002.

\bibitem[DrVa]{DrVa}
{\sc T.A. Driscoll} and {\sc S.A. Vavasis:}
{\sl Numerical conformal mapping 
using cross-ratios and Delaunay triangulation.}  SIAM J. Sci. Comput.  
19 (1998),  no. 6, 1783--1803 (electronic). 

\bibitem[DuVu]{DV}
{\sc V.~Dubinin} and {\sc M. Vuorinen,\,}
{\sl On conformal moduli of polygonal quadrilaterals}.
Helsinki preprint 417 August 2005.

\bibitem[Gai]{Gai2} {\sc D. Gaier,\,}
{\sl Conformal modules and their computation, }
in `Computational Methods and Function Theory' (CMFT'94),
R.M.Ali {\it et al.} eds., pp.159-171. World Scientific, 1995.


\bibitem[HVV]{HVV}
{\sc V. Heikkala, M.\ K. Vamanamurthy} and {\sc M. Vuorinen,\,}
{\sl Generalized elliptic integrals},
Helsinki preprint 404, 2004.

\bibitem[Hen]{Hen} {\sc P. Henrici,\,}
{\sl Applied and Computational Complex Analysis, vol. III,}
Wiley, Interscience, 1986.

\bibitem[Hu]{Hu} {\sc C. Hu,\,}
{\sl A software package for computing Schwarz-Christoffel conformal
transformation for doubly connected polygonal regions},
ACM Transactions of Math. Software 24 (1998).

\bibitem[IT]{IvTr} {\sc V.I. Ivanov} and {\sc M.K. Trubetskov,\,}
{\sl Handbook of Conformal Mapping with Computer-aided Visualization.}
CRC Press, 1995.

\bibitem[K\"uh]{Kuhnau:2005}
{\sc R. K\"uhnau,\,}
{\sl The conformal module of quadrilaterals and of rings,} In:
{ Handbook of Complex Analysis: Geometric Function Theory}, (ed.
by  R. K\"uhnau) Vol. 2,
North Holland/Elsevier, Amsterdam, 99-129, 2005.

\bibitem[LV]{Lehto:1973}
{\sc O. Lehto} and {\sc K.\ I. Virtanen,\,}
{\sl Quasiconformal mappings in the plane},
Springer, Berlin, 1973.


\bibitem[Opf]{Opf} {\sc G. Opfer,\,}
{\sl Die Bestimmung des Moduls zweifach zusammenh\"angender Gebiete
mit Hilfe von Differenzenverfahren}, Arch. Rat. Mech. Anal. 32 (1969),
281-297.

\bibitem[Pap]{Pap}
{\sc N. Papamichael,\,}
{\sl Dieter Gaier's contributions to numerical conformal mapping},
Comput. Methods Funct. Theory 3 (2003),  no. 1-2, 1-53.

\bibitem[PKo]{PaKo} {\sc N. Papamichael} and {\sc C.A. Kokkinos,\,}
{\sl The use of singular functions for the approximate conformal mapping
of doubly-connected domains,}
SIAM J. Sci. Stat. Comp. 5 (1984), 684-700.



\bibitem[PS]{PS}
{\sc N. Papamichael and N. S. Stylianopoulos:}
{\it The asymptotic behavior of conformal modules of quadrilaterals 
with applications to the estimation of resistance values,  }
Constr. Approx.  15  (1999),  no. 1, 109--134.



\bibitem[PWa]{PaWa} {\sc N. Papamichael} and {\sc M.K. Warby,\,}
{\sl Pole-type singularities and the numerical conformal mapping of
doubly-connected domains}, J. Comp. Appl. Math. 10 (1984), 93-106.


\bibitem[RV]{RV}
{\sc A. Rasila} and {\sc M. Vuorinen:}
Web-based capacity computation
{\it Work in progress}

\bibitem[Vla]{Vla} {\sc V.I. Vlasov,\,}
{\sl Multipole method for solving some boundary value problems in
complex-shaped domains,}
Zeitschr. Angew. Math. Mech. 76, suppl. 1 (1996), 279-282.

\bibitem[Weg]{Weg05}
{\sc R. Wegmann,\,}
{\sl Methods for numerical conformal mapping,  }
Handbook of complex analysis: geometric function theory.
Vol. 2,  (ed. by  R. K\"uhnau), Elsevier, Amsterdam,  351--477,  2005.

\bibitem[Wei1]{Wei1} {\sc J. Weisel,\,}
{\sl L\"osung singul\"aren Variationsprobleme durch die Verfahren von Ritz
und Galerkin mit finiten Elementen}, Anwendungen in der konformen
Abbildung,
Mitt. Math. Sem. Giessen 138 (1979), 1-150.

\bibitem[Wei2]{Wei2} {\sc J. Weisel,\,}
{\sl Approximation quasikonformer Abbildungen mehrfach
zusammenh\"angender Gebiete durch finite Elemente,}
J. Appl. Math 32 (1981), 34-44.

\end{thebibliography}
\end{document}